\documentclass[12pt]{article}
\usepackage[T2A]{fontenc}
\usepackage[cp1251]{inputenc}
\usepackage[russian]{babel}
\usepackage{amssymb,amsmath}
\usepackage{hyperref}
\usepackage{graphicx}
\usepackage{ccaption}
\usepackage{subfigure}
\captiondelim{. }
\usepackage{amsmath}
\usepackage{amssymb}
\usepackage{amsxtra}
\usepackage{amsthm}
\usepackage{latexsym}
\usepackage{array}
\usepackage{multirow}
\usepackage{hhline}
\usepackage{color}
\theoremstyle{plain}
\newtheorem{theorem}{Теорема}

\begin{document}

\noindent\hspace{0.55\textwidth}\parbox[t]{0.55\textwidth}{%
Дополнение к результатам Ф. Шарова\\
Рябов Павел Павлович\\[2ex]
}%
\footnotetext{Автор частично поддержан грантом Президента РФ МК-6137.2016.1}

В данной работе рассматривается вопрос о замощении прямоугольника меньшими прямоугольниками, отношения сторон которых принадлежат заданному набору квадратных иррациональностей. Эта проблема ранее была изучена Ф. Шаровым в работе~\cite{Sharov}. Сформулируем её основной результат.
\begin{theorem}[Ф. Шаров, 2016]
\label{theor01}
Пусть $x_1=a_1+b_1\sqrt{p},...,x_n=a_n+b_n\sqrt{p}$ --- такие числа, что $x_i>0, a_i,b_i,p\in{\mathbb Q}$, где $1\le i\le n$ и $\sqrt{p}\notin{\mathbb Q}$. Тогда:

1) если существуют такие числа $i$ и $j$, что $1\le i, j\le n$ и $(a_i-b_i\sqrt{p})(a_j-b_j\sqrt{p})<0$, то прямоугольник с отношением сторон $z$ можно разрезать на прямоугольники с отношениями сторон $x_1,...,x_n$ тогда и только тогда, когда $$z\in{\{e+f\sqrt{p}>0| e,f\in{\mathbb Q}}\};$$

2) если для всех $i$, где $1\le i \le n$, выполнено неравенство $a_i-b_i\sqrt{p}$>0, то прямоугольник с отношением сторон $z$ можно разрезать на прямоугольники с отношениями сторон $x_1,...,x_n$ тогда и только тогда, когда $$z\in \{e+f\sqrt{p}>0| e,f\in{\mathbb Q}, e>0, \frac{|f|}{e}\le \max\limits_{i}\frac{|b_i|}{a_i}\};$$

3) если для всех $i$, где $1\le i\le n$ выполнено неравенство $a_i-b_i\sqrt{p}<0$, то прямоугольник с отношением сторон $z$ можно разрезать на прямоугольники с отношениями сторон $x_1,...,x_n$ тогда и только тогда, когда $$z\in \{e+f\sqrt{p}>0| e,f\in{\mathbb{Q}, f>0, \frac{|e|}{f}\le\max\limits_i\frac{|a_i|}{b_i}}\}.$$
\end{theorem}
В статье~\cite{Sharov} доказано существование требуемых замощений, однако нет их явной конструкции. Здесь мы приведём такую конструкцию и получим близкий новый результат.

\section{Простое доказательство части ``тогда'' пункта 3) теоремы~\ref{theor01}}

Без ограничения общности положим $\max\limits_i\frac{|a_i|}{b_i}=\frac{|a_1|}{b_1}$ и $a_1, b_1, e, f\geq 0$. Действительно, из неравенств $a_1-b_1\sqrt{p}<0$, $a_1+b_1\sqrt{p}>0$ следует, что $b_1>0$. Если, например, $a_1<0$, тогда $\frac{1}{a_1+b_1\sqrt{p}}=\frac{-a_1}{b_1^2p-a_1^2}+\frac{b_1}{b_1^2p-a_1^2}\sqrt{p}=a_1'+b_1'\sqrt{p}$ --- отношение сторон того же прямоугольника, и уже $a_1', b_1'\ge 0$.

\begin{figure}[!ht]
\centering
\includegraphics[width=15cm,height=15cm, keepaspectratio]{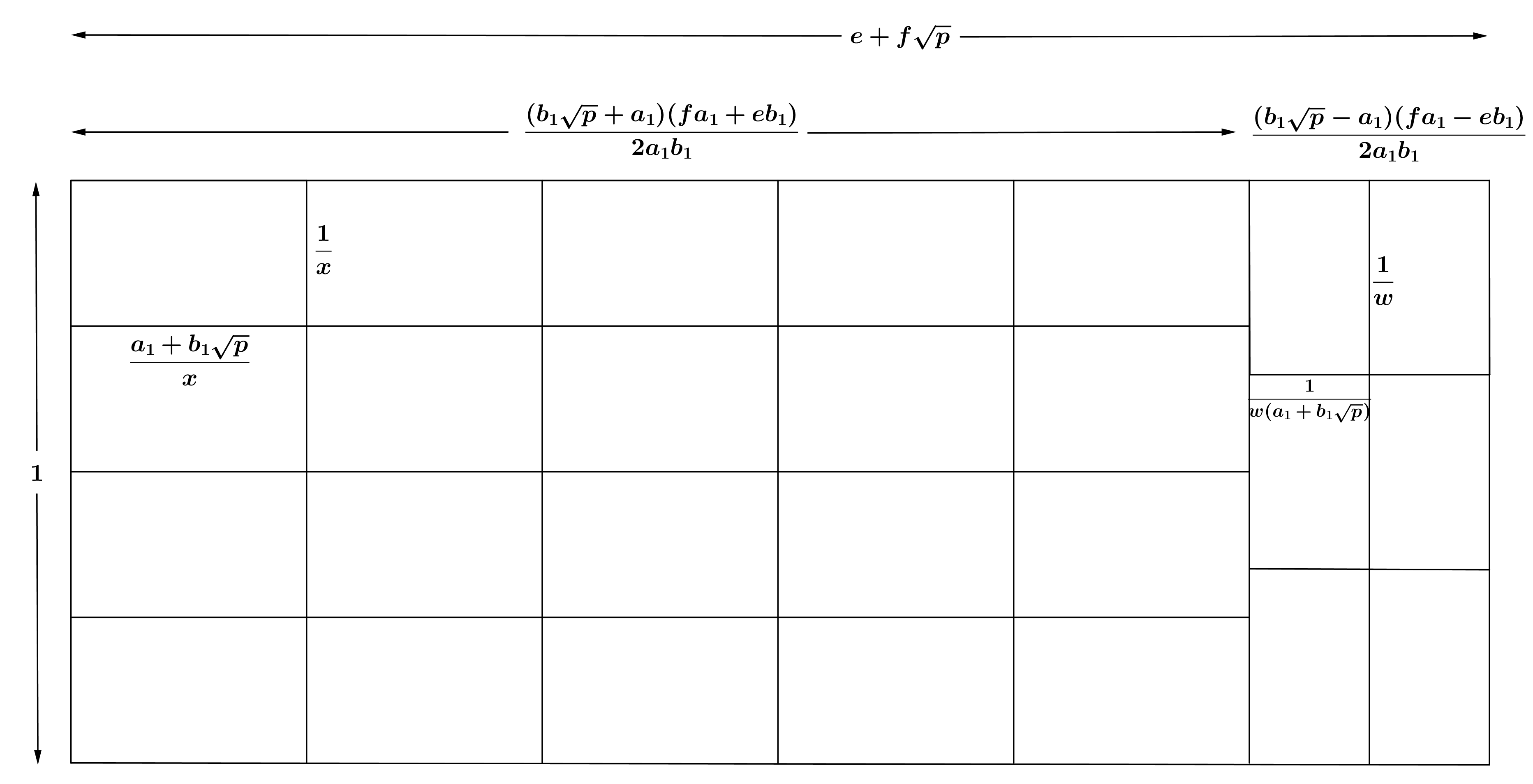}
\caption{Замощение в теореме~\ref{theor01} для $a_1=\frac{1}{2}, b_1=\frac{1}{2}, e=\frac{7}{24}, f=\frac{23}{24}, p=5$.}\label{figa01}
\end{figure}

Теперь рисунок~\ref{figa01} даёт явную конструкцию замощения прямоугольника $1\times(e+f\sqrt{p})$ прямоугольниками с отношением сторон $a_1+b_1\sqrt{p}$. Замощаемый прямоугольник разрезается на два прямоугольника $1\times\frac{(b_1\sqrt{p}+a_1)(fa_1+eb_1)}{2a_1b_1}$ и $1\times\frac{(b_1\sqrt{p}-a_1)(fa_1-eb_1)}{2a_1b_1}$. Положим $\frac{fa_1+eb_1}{2a_1b_1}=\frac{y}{x}$, а $\frac{(p{b_1}^2-{a_1}^2)(fa_1-eb_1)}{2a_1b_1}=\frac{u}{w}$, где $x, y, u, w\in\mathbb{Z}_+$. Прямоугольник $1\times\frac{(a_1+b_1\sqrt{p})(fa_1+eb_1)}{2a_1b_1}$ замостим $xy$ равными прямоугольниками $\frac{1}{x}\times\frac{a_1+b_1\sqrt{p}}{x}$, разделив его на $x$ равных частей по вертикали и на $y$ равных частей по горизонтали. Прямоугольник $1\times\frac{(b_1\sqrt{p}-a_1)(fa_1-eb_1)}{2a_1b_1}$ аналогично замостим $uw$ прямоугольниками $\frac{1}{w}\times\frac{1}{w(a_1+b_1\sqrt{p})}$.

Замощение в пунктах 1) и 2) теоремы~\ref{theor01} строится похожим образом.

\section{Новый результат}

Хотя в рассматриваемой задаче речь идёт о \emph{наборе} квадратичных иррациональностей, фактически в замощении используется только одна из них (см. $\S 1$). А что если потребовать, чтобы для каждого числа из набора в замощении нашелся хотя бы один прямоугольник с таким отношением сторон? Назовём разрезание прямоугольника на прямоугольники с отношениями сторон $x_1,..., x_n$ \emph{разнообразным}, если для каждого $i=1,...,n$ в замощении \emph{найдётся хотя бы один прямоугольник с отношением сторон} $x_i$. Мы покажем, что тогда выполняется следующая теорема (по сравнению с теоремой~\ref{theor01} изменено описание \emph{набора}, добавлено слово \emph{разнообразно} и в конце пунктов 2) и 3) неравенство заменено на строгое):
\begin{theorem}
\label{theor02}
Пусть $x_1=a_1+b_1\sqrt{p},...,x_n=a_n+b_n\sqrt{p}$ -- такие числа, что $\frac{x_i}{x_j}, x_i x_j\notin\mathbb{Q}$ при $i\neq j$, $x_i>0, a_i,b_i,p\in{\mathbb Q}$, где $1\le i\le n$ и $n\geq 2$. Тогда:

1) если существуют такие числа $i$ и $j$, что $1\le i, j\le n$ и $(a_i-b_i\sqrt{p})(a_j-b_j\sqrt{p})<0$, то прямоугольник с отношением сторон $z$ можно \emph{разнообразно} разрезать на прямоугольники с отношениями сторон $x_1,...,x_n$ тогда и только тогда, когда $$z\in{\{e+f\sqrt{p}>0| e,f\in{\mathbb Q}}\};$$

2) если для всех $i$, где $1\le i \le n$, выполнено неравенство $a_i-b_i\sqrt{p}$>0, то прямоугольник с отношением сторон $z$ можно \emph{разнообразно} разрезать на прямоугольники с отношениями сторон $x_1,...,x_n$ тогда и только тогда, когда $$z\in \{e+f\sqrt{p}>0| e,f\in{\mathbb Q}, e>0, \frac{|f|}{e}<\max\limits_{i}\frac{|b_i|}{a_i}\};$$

3) если для всех $i$, где $1\le i\le n$, выполнено неравенство $a_i-b_i\sqrt{p}<0$, то прямоугольник с отношением сторон $z$ можно \emph{разнообразно} разрезать на прямоугольники с отношениями сторон $x_1,...,x_n$ тогда и только тогда, когда $$z\in \{e+f\sqrt{p}>0| e,f\in{\mathbb{Q}, f>0, \frac{|e|}{f}<\max\limits_i\frac{|a_i|}{b_i}}\}.$$
\end{theorem}
\section{Доказательство теоремы~\ref{theor02}}

\begin{proof}[Доказательство пункта 1) теоремы~\ref{theor02}]

Необходимость следует из теоремы~\ref{theor01}. Докажем достаточность.
Разделим прямоугольник $1\times z$ на $n$ прямоугольников $A_1,...,A_n$ размера $\frac{1}{n}\times z$. К меньшей стороне прямоугольника $A_k$ изнутри приставим большей стороной прямоугольник с отношением сторон $x_k$ так, чтобы эти стороны полностью совпали. Оставшуюся часть прямоугольника $A_k$ замостим прямоугольниками с отношениями сторон $x_i$ и $x_j$, по теореме~\ref{theor01} это возможно. Так, для каждого $k$ в замощении найдется прямоугольник с отношением сторон $x_k$.
\end{proof}

В теореме~\ref{theor02} доказательство пунктов 2) и 3) аналогичны. Поэтому в данной работе будет приведено только доказательство пункта 2).

\begin{figure}[!ht]
\centering
\includegraphics[width=15cm,height=15cm, keepaspectratio]{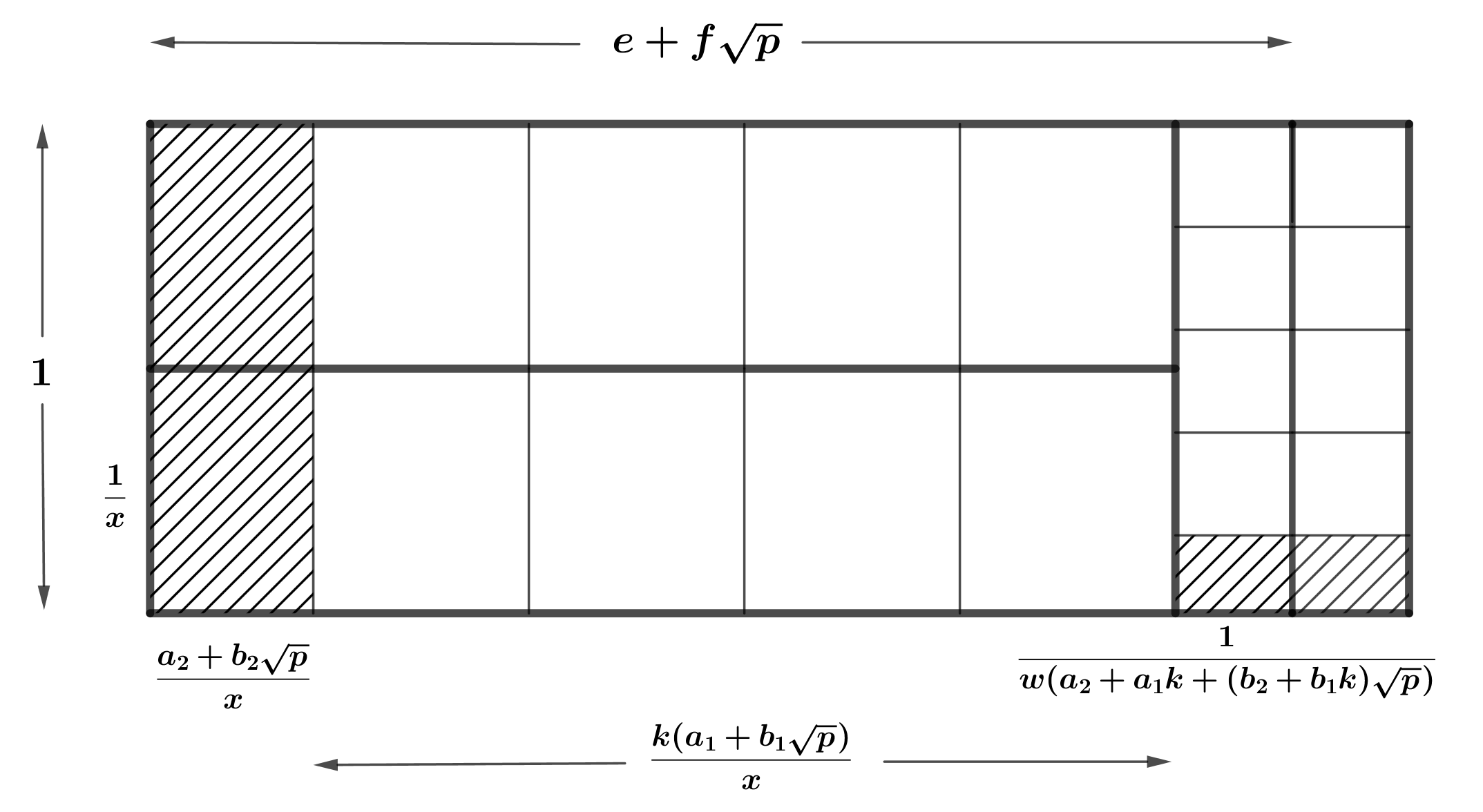}
\caption{}\label{figa02}
\end{figure}

\begin{proof}[Доказательство существования замощения в пункте 2) теоремы~\ref{theor02}]

Пусть $z=e+f\sqrt{p}$, где $e, f\in\mathbb{Q}$. Без ограничения общности $\frac{|f|}{e}<\frac{|b_1|}{a_1}$. Существует такое достаточно большое $k\in{\mathbb N}$, что $ka_1+a_2+...+a_n-(kb_1+b_2+...+b_n)\sqrt{p}>0$ и $\frac{|f|}{e}<\frac{|kb_1+b_2+...+b_n|}{ka_1+a_2+...+a_n}$, так как $k(a_1-b_1\sqrt{p})>0$ и $\lim\limits_{k\to +\infty} \frac{|kb_1+b_2+...+b_n|}{ka_1+a_2+...a_n}=\frac{|b_1|}{a_1}>\frac{|f|}{e}$. Прямоугольник с отношением сторон $r=ka_1+a_2+...+a_n+(kb_1+b_2+...+b_n)\sqrt{p}$ состоит из $k$ прямоугольников с отношением сторон $x_1$ и из прямоугольников с отношением сторон $x_2,...,x_n$. Прямоугольник $1\times z$ замостим прямоугольниками с отношением сторон $r$, по теореме~\ref{theor01} это возможно (а $\S 1$ даёт явную конструкцию замощения). На рисунке~\ref{figa02} представлено замощение для $n=2$ и $a_1, a_2, b_1, b_2, e, f\geq 0, \frac{b_1}{a_1}>\frac{f}{e}$, а числа $x, y, u, w\in\mathbb{Z}_+$ определяются равенствами $\frac{f(a_1k+a_2)+e(b_1k+b_2)}{2(a_1k+a_2)(b_1k+b_2)}=\frac{y}{x}$ и $\frac{((a_1k+a_2)^2-p(b_1k+b_2)^2)(e(b_1k+b_2)-f(a_1k+a_2))}{2(a_1k+a_2)(b_1k+b_2)}=\frac{u}{w}.$
\end{proof}

\begin{proof}[Доказательство невозможности замощения в пункте 2) теоремы~\ref{theor02}]

По теореме~\ref{theor01}, если прямоугольник $1\times z$ можно замостить прямоугольниками с отношениями сторон $x_1,...,x_n$, то $z=e+f\sqrt{p}$, для некоторых $e,f\in\mathbb{Q}$ таких, что $\max\limits_{i}\frac{|b_i|}{a_i}\geq\frac{|f|}{e}$. Остаётся доказать, что если $\max\limits_{i}\frac{|b_i|}{a_i}=\frac{|f|}{e}$, то прямоугольник $1\times z$ невозможно \emph{разнообразно} разрезать на прямоугольники с отношениями сторон $x_1,...,x_n$.

Предположим, что такое \emph{разнообразное} разрезание существует. Без ограничения общности будем считать $a_1, b_1, e, f$ положительными, аналогично \S 1,  $\frac{|b_2|}{a_2}<\frac{b_1}{a_1}=\frac{f}{e}$ ($\frac{|b_2|}{a_2}\neq\frac{b_1}{a_1}$ так как $\frac{x_i}{x_j}, x_i x_j\notin\mathbb{Q}$).
Измельчим это замощение следующим образом. Для каждого $i=2,...,n$ замостим каждый прямоугольник разбиения с отношением сторон $x_i$ прямоугольниками с отношением сторон $x_1$. Это возможно по теореме~\ref{theor01}, причём в явной конструкции такого замощения в $\S 1$, см. рисунок~\ref{figa01}, найдутся негомотетичные прямоугольники.

Дальнейшее доказательство существенно использует понятия и факты на страницах 209-210 статьи~\cite{Sharov}.

Стороны всех прямоугольников измельчения представляются в виде $\alpha+\beta\sqrt{p}$, где $\alpha,\beta\in\mathbb{Q}$ (см. третий абзац стр. 210).

Рассмотрим ``площади'' (определение 5 на стр. 209) прямоугольников измельчения. ``Площадь'' прямоугольника $\alpha+\beta\sqrt{p}\times x_1(\alpha+\beta\sqrt{p})$, где $\alpha, \beta\in\mathbb{Q}$, равна $\beta^2\left(\frac{2fa_1^3}{b_1^2}-pfa_1-peb_1\right)=\beta^2\left(\frac{2fa_1(a_1^2-pb_1^2)}{b_1^2}\right)\geq 0$ (стр. 210, формула для $S_M$).

Если хотя бы у одного прямоугольника $\beta\neq 0$, то общая сумма ``площадей'' прямоугольников измельчения строго больше $0$. Однако ``площадь'' замощаемого прямоугольника равна $0$ (стр. 210, формула для $S_{\text{Б}}$), противоречие.

Остаётся рассмотреть случай, когда все прямоугольники измельчения имеют вид $\alpha\times x_1\alpha$, где $\alpha\in\mathbb{Q}$. Будем считать, что у измельчаемого прямоугольника сторона $1$ вертикальна. Рассмотрим прямоугольник $A_1$, у которого сторона $\alpha x_1$ вертикальна (такой найдется, так как в измельчении есть негомотетичные прямоугольники). Выберем прямоугольник $A_2$ над $A_1$, имеющий с $A_1$ общий отрезок, для $A_2$ так же выберем прямоугольник $A_3$ над $A_2$, имеющий с $A_2$ общий отрезок, и т.д. Аналогично выберем прямоугольники под $A_1$. Сумма всех их вертикальных сторон будет равняться вертикальной стороне измельчаемого прямоугольника, то есть $1$. Эта же сумма будет равна $x_1w+u$, для некоторых $w, u\in\mathbb{Q}, u\geq 0, w>0$. Получим равенство $1=w(a_1+b_1\sqrt{p})+u$. Противоречие, так как $w, b_1\neq0$ по строгому неравенству $\frac{|b_2|}{a_2}<\frac{b_1}{a_1}$.
\end{proof}
\section*{Благодарности}
Этой работы бы не было без наставлений М. Б. Скопенкова. Автор благодарен Ф. А. Шарову за ценные замечания.

\end{document}